\documentclass[12pt]{elsarticle}
\usepackage{amsfonts}
\usepackage{amssymb}
\usepackage{mathrsfs}
\usepackage{amsmath}
\usepackage{hyperref}

\def\N{\mathbb{N}}

\def\cD{\mathcal D}
\def\cP{\mathcal P}
\def\cH{\mathcal H}
\def\hc{\hat c}

\def\pr{^\prime }
\def\tr{^{\mathsf {T}}}
\def\eop{\unskip\nobreak\hfil\penalty50\hskip2em\hbox{}\nobreak
\hfill\mbox{$\Box $}\par}
\newtheorem{theorem}{Theorem}[section]

\newtheorem{lemma}{Lemma}[section]
\begin{document}
\begin{frontmatter}
	\title{A unified construction of all the hypergeometric and basic hypergeometric families of orthogonal polynomial sequences}
	
\author{Luis Verde-Star}
 \address{
Department of Mathematics, Universidad Aut\'onoma Metropolitana, Iztapalapa,
Apartado 55-534, Mexico City 09340,
 Mexico }
\ead{verde@xanum.uam.mx}
	\begin{abstract}
	We construct a set $\cH$ of orthogonal polynomial sequences that contains all the families in the Askey scheme and the $q$-Askey scheme. The polynomial sequences in $\cH$ are solutions of a generalized first-order difference equation which is determined by three linearly recurrent sequences of numbers. Two of these sequences are solutions of the difference equation $s_{k+3}= z\, (s_{k+2} - s_{k+1}) + s_k$, where $z$ is a complex parameter, and the other sequence satisfies a related difference equation of order five. 

 We obtain explicit expressions for the coefficients of the orthogonal polynomials and for the generalized moments with respect to a basis of Newton type of the space of polynomials. We also obtain explicit formulas for the coefficients of the three-term recurrence relation satisfied by the polynomial sequences in $\cH$.

	The set $\cH$ contains all the 15 families in the Askey scheme of hypergeometric orthogonal polynomials \cite[p. 183]{Hyp} and all the 29 families of basic hypergeometric orthogonal polynomial sequences in the $q$-Askey scheme \cite[p. 413]{Hyp}. Each of these families is obtained by direct substitution of appropriate values for the parameters in our general formulas. The only cases that require some limits are the Hermite and continuous $q$-Hermite polynomials. We present the values of the parameters for some of the families.


\end{abstract}
\end{frontmatter}

\section{Introduction}
Among the families of orthogonal polynomial sequences the hypergeometric and basic hypergeometric families are certainly some of the most important and have been extensively studied for a long time. See \cite{GasRa}, \cite{Ism}, and \cite{Hyp}. 
In the present paper, we present a construction of a class $\cH$ of orthogonal polynomial sequences that includes all the hypergeometric and basic hypergeometric families. 
Our construction depends on three sequences of numbers, one of them determines a Newton basis for the space of polynomials, and the other two sequences determine a linear operator on the space on polynomials defined by its action on the Newton basis. The eigenfunctions of the operator form a polynomial sequence. If the three sequences satisfy certain difference equations and their initial terms are related in a suitable way then the sequence of polynomial eigenfunctions satisfies a three-term recurrence relation and it is orthogonal with respect to a positive measure, or with respect to a moments functional.

 We obtain explicit expressions for the coefficients of the orthogonal polynomials and for the generalized moments with respect to the Newton basis. We also obtain explicit formulas for the coefficients of the three-term recurrence relation satisfied by the polynomial sequences of eigenfunctions in $\cH$. 
 We have verified that $\cH$ contains all the 15 families in the Askey scheme and all the 29 families in the $q$-Askey scheme. The coefficients of the normalized three-term recurrence relation satisfied by each family can be obtained by substitution of suitable values of the parameters in our general formulas for the coefficients of the recurrence relations satisfied by the elements of $\cH$. The only cases that require taking a limit of one of the parameters are the Hermite and the continuous $q$-Hermite polynomials.

We present some examples of cases in which the generalized difference equation of first order becomes the usual second order differential or difference equation or $q$-difference equation used to characterize some families of orthogonal polynomial sequences.

Since our results produce a uniform parametrization of all the hypergeometric and basic hypergeometric orthogonal polynomial sequences they may be useful to study the structure of the space of such sequences. See \cite{Koor}.

In \cite{Ger} Geronimus initiated the study of orthogonal polynomial sequences expressed in terms of a Newton basis. Al-Salam and Verma \cite{AlSV} considered the case in which the nodes of the Newton polynomials are of the form $x_k= a + b q^{-k}$. Recently, Vinet and Zhedanov \cite{VZ} studied the cases where the sequences of eigenvalues correspond to the classical grids and obtained results that are very similar to some of our results, but they used a different approach.  

There are numerous recent papers that deal with recurrence relations, moments, characterization theorems, difference operators, and orthogonality on quadratic lattices of polynomial sequences related to the ones we study in this paper. Some of them use infinite matrices. See for example, \cite{CosMar}, \cite{KS}, \cite{MBP}, \cite{Mas}, \cite{quad}, \cite{Mom}, and \cite{Tch}.

In the following section we present the main ideas used in our development and at the end of the section we describe the further contents of the paper.  

\section{Main steps in the construction of $\cH$.}
Here we present a brief description of the main steps in our construction of the set $\cH$ of orthogonal polynomial sequences. A more detailed account and the proofs of the statements in this section are presented in the subsequent sections.

Let $h_k$ and $x_k$ be sequences of numbers that satisfy the recurrence relation  
 $$ s_{k+3}=z (s_{k+2}-s_{k+1}) + s_k, \qquad k \ge 0, \eqno(2.1) $$
 where $z$ is a complex parameter, 
 and let $g_k$ be a sequence with $g_0=0$ that satisfies the recurrence relation  
 $$s_{k+5}=(z^2-z-1) (s_{k+4}-s_{k+1}) -(z-1) (z^2-z-1) (s_{k+3}-s_{k+2}) + s_k, \qquad k \ge 0. \eqno(2.2)$$
 This recurrence relation is satisfied by the termwise product of any pair of solutions of (2.1).

Let $\{v_k \}$ be the Newton basis associated with the sequence $\{x_k \}$, defined by 
 $v_0(t)=1$, and
$$v_k(t)= (t-x_0) (t-x_1) (t-x_2) \cdots (t-x_{k-1}), \qquad k\ge 1. $$

Let $\cD$ be the linear operator on the space of polynomials determined by $\cD v_k = h_k v_k + g_k v_{k-1}$, since $g_0=0 $ we see that $\cD t^n = h_n t^n +$ polynomial of lower degree. Define the polynomial $u_n$ as a monic polynomial of degree $n$ which is an eigenfunction of $\cD$ with eigenvalue $h_n$. That is 
$$\cD u_k = h_k u_k, \qquad k \ge 0. \eqno(2.3).$$

The operator $\cD$ is a generalized difference operator which in concrete examples becomes a second order differential or difference or $q$-difference operator on a linear or quadratic lattice. We will show that 
 $$u_n(t) = \sum_{k=0}^n c_{n,k} v_k(t), \qquad n \ge 0,\eqno(2.4) $$
 holds with $ c_{n,k}$ given by
$$c_{n,k}=\prod_{j=k}^{n-1} \frac{g_{j+1}}{h_n -h_j}, \qquad 0 \le k \le n-1, \eqno(2.5) $$
and $c_{n,n}=1$ for $n \ge 0$.
 In concrete examples this becomes a ($q$-)hypergeometric representation for the polynomials $u_n (t)$.

If the three sequences are connected in a suitable way then the polynomial sequence $\{u_n(t)\}$ satisfies a three-term recurrence relation of the form 
$$u_{n+1}(t)= (t-\beta_n) u_n(t) - \alpha_n u_{n-1}, \qquad n \ge 1, \eqno(2.6)$$
and also satisfies the generalized difference equation (2.3).
If all the $\alpha_n$ are positive and the $\beta_k$ are real then the sequence $\{u_n\}$ is orthogonal with respect to a positive measure, and if all the $\alpha_n$ are nonzero then $\{u_n\}$ is orthogonal with respect to a not necessarily positive definite moments functional.

 We will obtain explicit expressions in terms of the $h_k, x_k, g_k$ for the coefficients $\alpha_n$ and $\beta_n$ of the three-term recurrence relation (2.6) in Section 5. 

The generalized moments of the polynomials $u_k(t)$ with respect to the basis $\{ v_k(t) \}$ are the entries in the 0-th column of the inverse of the matrix $[c_{n,k}]$, and are given by
$$ m_n =\prod_{k=1}^n \frac{g_k}{h_0-h_k}, \qquad n \ge 1, \eqno(2.7)$$
and $m_0=1$. Note that they satisfy a recurrence relation of order one. If $g_{N+1}=0$ for some $N$ then $m_n=0$ for $n>N$. In such case the polynomials $u_k(t)$ are defined only for $0 \le k \le N$.
The standard moments with respect to the basis of monomials $\{t^n\}$ are easily obtained using a change of basis matrix.

The class $\cH$ is the set of all the orthogonal polynomial sequences $\{u_k(t)\}$ obtained by the procedure described above. The polynomial sequences in $\cH$ are determined by $z$ and the initial values of the sequences $\{h_k\}$, $\{x_k\}$, and $\{g_k\}$. 
 The sequences $\{u_k(t)\}$ are well defined if the $h_k$ are pairwise distinct and $g_k \ne 0$ for $k \ge 1.$ 

The roots of the characteristic polynomial of the difference equation (2.1) are $1, q$, and $ q^{-1}$, where $q$ is a nonzero complex number and $z=1+q+q^{-1}$. We classify the elements of $\cH$ by considering the possible multiplicities of the roots. There are three cases:
\begin{enumerate}
	\item \  The 3 roots are distinct, that is, $q\ne 1$ and $q\ne -1$.

	\item \  $q=1$, that is, 1 is a root with multiplicity 3.

	\item \  $q=-1$, that is, $-1$ is a double root and 1 is a simple root.
\end{enumerate}

Let $\cH_q$, $\cH_1$, and $\cH_{-1}$ be the subsets of $\cH$ that correspond to the first, second, and third cases respectively. The elements of $\cH_q$ are the basic hypergeometric polynomial sequences in the $q$-Askey scheme. 
 $\cH_1$ contains the hypergeometric orthogonal polynomial sequences in the Askey scheme. 
 The set $\cH_{-1}$ contains some polynomial sequences that are not as well-known as those in $\cH_q$ and $\cH_1$. Some particular elements of this set have been studied recently in \cite{Vinet1}, \cite{Vinet2}, and \cite{Vinet3}.

For the class $\cH_q$, we have $q \ne 1$ and $q \ne -1$, and the sequences $\{h_k\}$, $\{x_k\}$, and $\{g_k\}$ can be expressed as 
$$h_k= a_0 + a_1 q^k + a_2 q^{-k}, \qquad k \ge 0, \eqno(2.8)$$
$$x_k= b_0 + b_1 q^k + b_2 q^{-k}, \qquad k \ge 0, \eqno(2.9)$$
and
$$g_k= d_0 + d_1 q^k + d_2 q^{-k} + d_3 q^{2k} + d_4 q^{-2k} , \qquad k \ge 0, \eqno(2.10)$$
where the relations
$$d_0= -( a_2 b_2 q + d_1 + d_2 + a_1 b_1 q^{-1}),\qquad d_3=a_1 b_1 q^{-1}, \qquad d_4= a_2 b_2 q, \qquad \eqno(2.11)$$
are required to obtain orthogonal polynomial sequences $\{u_k(t)\}$.
Using the explicit formulas of sections 7 and 8 for the coefficients of the three-term recurrence relation (2.6) in terms of the parameters $a_0,a_1,a_2,b_0,b_1,b_2,d_1,d_2$ we can obtain the normalized recurrence relation of each family in the $q$-Askey scheme by giving suitable values to the parameters. For example, for the Askey-Wilson family we have
$$ a_1= a b c d q^{-1} a_2, \qquad b_0=0, \qquad b_1=a/2, \qquad b_2=a^{-1}/2, $$
$$ d_1=- a ( abcd +q (b c + b d + c d)) q^{-2} a_2/2, \qquad d_2=- ( (b + c + d) + q a^{-1} ) a_2/2, \eqno(2.12)$$
where $a_2$ is an arbitrary nonzero number and $a,b,c,d$ are the parameters in \cite[eq. 14.1.5]{Hyp}.

In Section 3 we present preliminary material and some definitions and notation. In Section 4 we find explicit expressions for the polynomial eigenfunctions $\{u_k\}$ of the operator $\cD$.
In Section 5 we find the matrix $L$ that represents the operator of multiplication by the independent variable with respect to the basis $\{u_k\}$.
In Section 6 we consider sequences $h_k$ and $x_k$, that satisfy (2.1), and $g_k$, that satisfies (2.2), and we find certain conditions on the initial terms of those sequences that make the matrix $L$ tridiagonal. We also state some recurrence relations satisfied by the entries of $L$.
In Section 7 we study the class $\cH_q$, which contains all the families in the $q$-Askey scheme, and we find values of the parameters that yield some of the families.
In Section 8 we consider the class $\cH_1$, which contains all the families in the Askey scheme, and we find values of the parameters that produce some of the families.
In Section 9 we look at the set $\cH_{-1}$, which corresponds to $q=-1$. Section 10 contains the proof of the theorem about the tridiagonality of the matrix $L$ stated in Section 6. 

\section{Preliminary material}
 We present in this section some preliminary material and introduce notation that will be used in the paper. A more detailed account of the matrix approach to polynomial sequences can be found in \cite{Mops} and \cite{Rec}. See also \cite{Aren} where some properties of doubly-infinite matrices are obtained. 

A polynomial sequence is a sequence of polynomials $p_0(t), p_1(t),p_2(t),\ldots $ with complex coefficients such that $p_n(t)$ has degree $n$ for $n \ge 0$. Every polynomial sequence is a basis for the complex vector space $\cP$ of all polynomials in one variable.

If $\{u_n\}$ and $\{v_n\}$ are two polynomial sequences then there exists a unique matrix $A=[a_{n,k}]$, where $(n,k) \in \N \times \N$, such that 
$$u_n(t) = \sum_{k=0}^n a_{n,k} v_k(t), \qquad n \ge 0. \eqno(3.1)$$
The infinite matrix $A$ is lower triangular and invertible. If all the polynomials $u_n$ and $v_n$ are monic then $a_{n,n}=1$ for $n\ge 0$.
If we consider a fixed polynomial sequence $v_n$ then every lower triangular invertible matrix $A$ determines another polynomial sequence by (3.1).
The $n$-th row of the matrix $A$ is the vector of coefficients of $u_n$ with respect to the basis $\{v_n\}$. Equation (3.1) is equivalent to the matrix equation
$$ [u_0,u_1,u_2,\ldots]\tr = A \ [v_0,v_1,v_2,\ldots ]\tr .\eqno(3.2)$$
We say that $A$ is the matrix of the sequence of polynomials $\{u_k(t): k \in \N\}$ with respect to the basis $\{v_k(t):k \in \N \}$.

Let $x_0, x_1,x_2,\ldots $ be a sequence of complex numbers and define the polynomials $v_0(t)=1$, and
$$v_k(t)= (t-x_0) (t-x_1) (t-x_2) \cdots (t-x_{k-1}), \qquad k\ge 1. \eqno(3.3)$$
It is clear that $\{v_n\}$ is a basis for the space of polynomials. It is called the {\sl Newton basis} associated with the sequence $x_n$. Let $V$ be the infinite matrix that satisfies
$$[v_0(t), v_1(t),v_2(t),\ldots ]\tr = V [1,t,t^2, \ldots]\tr. \eqno(3.4)$$
The entries in the $n$-row of $V$ are the coefficients of $v_n(t)$ with respect to the basis of monomials and therefore they are elementary symmetric functions of $x_0,x_1,x_2, \ldots, x_{n-1}$, with the appropriate signs. Therefore the entries in $V^{-1}$ are complete homogeneous symmetric functions of the $x_k$.

The dual basis of the Newton basis $\{v_n\}$ is the sequence of divided difference functionals $\Delta[x_0,x_1,x_{n-1}]$, which give us the coefficients in the representation of any polynomial in terms of the Newton basis. The basic theory of divided differences can be found in \cite{CdB} and \cite{ddci}.

Let us note that we consider the vector of coefficients of a polynomial $p$ with respect to a basis $\{v_k\}$ as a {\em row vector}, and thus the vectors of coefficients of a polynomial sequence with respect to a basis $\{v_k\}$ form a lower triangular infinite matrix. This convention is not standard and has as a consequence that applying a linear operator $\tau$ represented by a matrix $T$ with respect to a basis $\{v_k\}$, to a polynomial sequence represented by a matrix $A$, corresponds with multiplication by $T$ on the right, that is $AT$. Therefore, if $\tau_1$ and $\tau_2$ are operators whose matrix representations with respect to a basis $\{v_k\}$ are $T_1$ and $T_2$ respectively, then the composition of operators $\tau_2 \tau_1$ is represented by the matrix $T_1 T_2$ with respect to the same basis. That is, the map from operators to matrices is an anti-homomorphism. 

Let $\tau$ be a linear operator on the space of polynomials and let $T$ be its matrix representation with respect to a basis $\{v_k(t):k \in \N \}$, and let $A$ is the matrix of the sequence of polynomials $\{u_k(t): k \in \N\}$ with respect to the basis $\{v_k(t):k \in \N \}$ then $A T$ is the matrix of the sequence $\{ \tau u_k(t):k \in \N\}$ with respect to the same basis.

We introduce next some infinite matrices that will be used in the rest of the paper.
Let
$$ S= \left[ \begin{matrix}  0\  & 0\  & 0\  & 0\  & \ldots \cr
             1\  & 0\  & 0\  & 0\  & \ldots \cr
                 0\  & 1 \  & 0 \  & 0\  & \ldots \cr
                         0\  & 0 & 1 \  & 0\  & \ldots \cr
            \vdots \  & \vdots \  & \vdots \ & \ddots \  & \ \ddots \ \end{matrix}\right], \qquad 
 S\tr = \left[ \begin{matrix}  0\  & 1\  & 0\  & 0\  & \ldots \cr
                         0\  & 0\  & 1\  & 0\  & \ldots \cr
                        0\  & 0\  & 0 \  & 1\  & \ldots \cr
                        0\  & 0 & 0\  & 0\  & \ddots \cr
 \vdots \  & \vdots \  & \vdots \ & \vdots \  & \ \ddots \ \end{matrix}\right].  \eqno(3.5) $$
The matrix $S$ is called the left shift and $S\tr$ is the right shift.

Any sequence of numbers $g_0, g_1,g_2,\ldots $ can be used to construct an infinite matrix $G$ defined by $G_{k,k}= g_k$ for $k \ge 0$ and $G_{j,k}=0$ if $j \ne k$. We say that $G$ is the diagonal matrix associated with the sequence $g_k$. 
A matrix of the form $G S$, where $G$ is diagonal and $G_{k,k}\ne 0$ for $k \ge 1$, is called {\em generalized difference matrix} of first order. If the linear operator represented by $G S$ with respect to the basis $\{v_k(t)\}$ is denoted by $\gamma$ then 
$\gamma v_k(t)= g_k v_{k-1}(t)$, for $k \ge 1$, and $\gamma v_0(t)=0$.
 A diagonal matrix $G$, acting by multiplication on the right-hand side, represents an operator of the form $\phi v_k(t)= g_k v_k(t)$, which is a rescaling of the basis and can be considered as a generalized difference operator of order zero.  

Define the matrix  
$$ D= \left[ \begin{matrix}  0\  & 0\  & 0\  & 0\  & \ldots \cr
	        1\  & 0\  & 0\  & 0\  & \ldots \cr
		           0\  & 2 \  & 0 \  & 0\  & \ldots \cr
		              0\  & 0 & 3 \  & 0\  & \ldots \cr
		      \vdots \  & \vdots \  & \vdots \ & \ddots \  & \ \ddots \ \end{matrix}\right].  \eqno(3.6)$$
With respect to the standard basis of monomials $D$ is the matrix representation of the usual differential operator, and $S\tr$ is the representation of multiplication by the variable $t$ of the polynomials.

If $\{u_k\}_{k\geqslant 0}$ is a sequence of monic orthogonal polynomials we write the corresponding three-term recurrence relation in the form
 $$\alpha_k u_{k-1}(t)+\beta_k u_k(t)+u_{k+1}(t)= t u_k(t), \qquad k \geqslant 1, \eqno(3.7) $$
 and define the sums $\sigma_k= \beta_0 +\beta_1 +\cdots +\beta_k,$ for $k \geqslant 0.$ Using the $\sigma_k$ is convenient because their explicit formulas are simpler than those for the $\beta_k$.

\section{The generalized difference equation}
We denote by $\cP$ the complex vector space of all polynomials in one variable.
Let $v_k(t)$ be a monic polynomial sequence in $\cP$. It is clear that $\{v_k: k\ge 0\}$ is a basis for the space $\cP$. For our purposes, a generalized difference operator of order one with respect to the basis $\{ v_k: k \ge 0\}$ is a linear operator $\gamma$ defined by $\gamma v_0=0$ and $\gamma v_k = g_k v_{k-1}$ for $k \ge 0$, where $g_k$ is a sequence of complex numbers.  
 A linear operator $\phi$ defined by $\phi v_k = h_k v_k$, where the $h_k$ are complex numbers, is called generalized difference operator of order zero.

Let $\gamma$ and $\phi$ be generalized difference operators with respect to the basis $\{v_k\}$, of orders one and zero, respectively, and let $\{u_n(t): n \ge 0\}$ be a monic polynomial sequence.
Let us consider the generalized difference equation
$$ \gamma u_n(t) + \phi u_n(t) = h_n u_n(t), \qquad n \ge 0. \eqno(4.1)$$
This is the difference-eigenvalue equation that we expressed as $\cD u_n = h_n u_n$ in section 2.
Note that if $h_k\ne 0$ for $k\ge 0$ then the operator $\cD$ is invertible.

\begin{theorem} If the sequence $h_k$ satisfies $ h_k \ne h_j$, for $k \ne j$, then the solution of (4.1) is the polynomial sequence 
	$$ u_n(t) = \sum_{k=0}^n c_{n,k} v_k(t), \qquad n \ge 0, \eqno(4.2) $$	
where the coefficients $c_{n,k}	$ are given by 
	$$ c_{n,k} =\prod_{j=k}^{n-1} \frac{ g_{j+1}}{ h_n - h_j}, \qquad 0 \le k \le n-1, \eqno(4.3)$$
	 and $c_{n,n}=1$ for $n\ge 0$. 
\end{theorem} 

{\it Proof:} Since $\{u_n(t): n \ge 0\}$ is a monic polynomial sequence it is clear that $u_n(t)$ can be written as in equation (4.2) for some matrix of coefficients $c_{n,k}$ with $c_{n,n}=1$ for $n\ge 0$. Then, by the definition of $\gamma$ we have
$$ \gamma u_n(t)= \sum_{k=1}^n c_{n,k} g_k v_{k-1}(t) = \sum_{k=0}^{n-1} c_{n,k+1} g_{k+1} v_k(t), $$
and the definition of $\phi$ gives
$$ \phi u_n(t)= \sum_{k=0}^n c_{n,k} h_k v_k(t). $$ 
Therefore the difference equation (4.1) becomes
$$\sum_{k=0}^{n-1} \left( c_{n,k+1} g_{k+1} + c_{n,k} h_k\right) v_k(t) +h_n v_n(t) = \sum_{k=0}^{n-1} c_{n,k} h_n v_k(t) + h_n v_n(t), $$  
and this gives
$$ \sum_{k=0}^{n-1} \left( c_{n,k+1} g_{k+1} + c_{n,k} (h_k - h_n) \right) v_k(t) = 0. $$
By the linear independence of the polynomial sequence $\{ v_k\}$ we obtain
 $$ c_{n,k+1} g_{k+1} = c_{n,k} (h_n - h_k), \qquad 0 \le k \le n-1, $$ 
 and since the numbers $h_j$ are pairwise distinct we can write the previous equation in the form  
$$ c_{n,k}= \frac{g_{k+1}}{h_n - h_k} \  c_{n,k+1}, \qquad  0 \le k \le n-1. $$ 
 This recurrence relation clearly gives us (4.3). \eop

 A result equivalent to (4.3) is obtained in \cite{VZ} using a different approach.
  
 \medskip
 Let $C$ be the matrix of coefficients $c_{n,k}$. It is an invertible lower triangular infinite matrix that satisfies
$$ C [v_0(t), v_1(t),v_2(t), \ldots]\tr = [u_0(t), u_1(t), u_2(t), \ldots ]\tr, \qquad n \ge 0. \eqno(4.4)$$
Let $V$ be the matrix whose $(n, k)$ entry is the Taylor coefficient of $t^k$ in the polynomial $v_n(t)$.
Then $ CV$ is the matrix of coefficients of the polynomials $u_k(t)$ with respect
to the basis of monomials $ \{t^k : k \ge 0\}$.

We define the polynomial sequence $w_k(t)$ as follows, $w_0(t)=1$ and
$$ w_k(t) = (t-h_0) (t-h_1) (t-h_2) \cdots (t- h_{k-1}), \qquad k\ge 1. \eqno(4.5)$$
It is the Newton basis associated with the sequence of nodes $h_k$. We will use the notation
$$w_{n,k}(t)=\frac{w_n(t)}{w_k(t)} = \prod_{j=k}^{n-1} (t-h_j), \qquad 0 \le k \le n. \eqno(4.6)$$

\begin{theorem}
	Let $C^{-1}=[\hc_{n,k}]$. Then	
	$$ \hc_{n,k}=\frac{ \prod_{j=k+1}^n g_j}{ w_{n+1,k}\pr (h_k) }=\prod_{j=k+1}^n \frac{g_j}{h_k - h_j}, \qquad 0 \le k \le n-1, \eqno(4.7)$$
	 and $\hc_{n,n}=1$ for $n \ge 0.$
\end{theorem}

{\it Proof:} Let us note that the denominator in (4.3) is equal to $ w_{n+1,k}\pr(h_n)$ if $ k < n$.
Let $k < n$ and let $ \hc_{n,j}$ be defined by (4.7). Then 
$$ \sum_{j=k}^n \hc_{n,j} c_{j,k} = \sum_{j=k}^n \frac{\left( \prod_{i=j+1}^n g_i \right) \left( \prod_{i=k+1}^j g_i \right)}{w_{n+1,j}\pr(h_j) \  w_{j+1,k}\pr(h_j)}= \left( \prod_{i=k+1}^n g_i \right) \sum_{j=k}^n \frac{1}{w_{n+1,k}\pr (h_j)} = 0, \eqno(4.8)$$ 
because the last sum is the divided difference of 1 with respect to $ h_k, h_{k+1},\ldots,h_n$, which are at least two nodes, since $k<n$. Alternatively we can use that the sum is equal to the sum of the residues of $1/w_{n+1,k}(t)$. The basic properties of divided differences can be found in \cite{CdB} or \cite{ddci}.  
 
 For $k=n$ we get $\hc_{n,n} c_{n,n} =1$ for $n \ge 0$. Therefore the matrix product $[\hc_{n,k}] [c_{n,k}]$ is equal to the infinite identity matrix and this completes the proof. \eop
 
The entries in the 0-th column of $C^{-1}$ are 
$$ \hc_{n,0} =\prod_{k=1}^n \frac{g_k}{h_0-h_k}, \qquad n \ge 1, \eqno(4.9)$$
and $\hc_{0,0}=1$. We denote them by $m_n=\hc_{n,0}$ for $ n \ge 0.$ They are especially important because if the polynomial sequence $u_k(t)$ is orthogonal with respect to some linear functional then $m_n$ is the generalized moment of $u_n$ with respect to the Newton basis $\{ v_k(t): k \ge 0 \}$. 
Note that the sequence $m_n$ satisfies a recurrence relation of order one.

\section{The operator of multiplication by the variable $t$}
Let $\{ v_k \}$ be the Newton basis associated with a sequence $x_0,x_1,x_2,\ldots$ as defined in (3.3), 
and let $V$ be the matrix whose $(n,k)$ entry is the Taylor coefficient of $t^k$ in the polynomial $v_n(t)$.
The entries of the inverse matrix $V^{-1}$ are the complete homogeneous symmetric polynomials of the nodes $x_j$, that is, 
$$ (V^{-1})_{n,k}= \sum x_{i_0} x_{i_1} \cdots x_{i_{n-k}}, \eqno(5.1)$$
where the sum runs over all vectors $ (i_0,i_1,\ldots, i_{n-k})$ with entries in $\{x_0, x_1, \ldots, x_k\}$. See \cite[p.21]{Mac}. 
Therefore $V$ and $V^{-1}$ are change of bases matrices that satisfy
$$ V [1, t, t^2, \ldots]\tr = [v_0(t), v_1(t), v_2(t),\ldots ]\tr , $$
and 
$$ V^{-1} [v_0(t), v_1(t), v_2(t),\ldots ]\tr = [1, t, t^2, \ldots]\tr . $$ 

In the previous section we saw that $C V$ is the matrix of coefficients of the polynomials $u_k(t)$ with respect to the basis of monomials $\{t^k : k \ge 0\}$, for any basis $\{v_k\}$, not necessarily of Newton type.

Since 
$$ t v_k(t) =( ( t - x_k) + x_k) v_k(t) = v_{k+1}(t) + x_k v_k(t), \qquad k \ge 0, $$
the matrix representation with respect to the basis $\{ v_k(t): k \ge 0\} $ of the map $ p(t) \rightarrow t p(t)$ on the space $\cP$ is $ S\tr + F$, where $F$ is the diagonal matrix whose $(k,k)$ entry is $x_k$, for $k \ge 0$.

Let $L$ be the matrix representation with respect to the basis $\{ u_k(t): k \ge 0\} $ of the operator of multiplication by $t$. Using equation (3.4) we obtain
$$ L = C ( S\tr + F ) C^{-1} . \eqno(5.2)$$

The explicit expressions (4.3) and (4.7) for the entries of $C$ and $C^{-1}$ give us  
$$ L_{n,k}=\left( \prod_{j=k+1}^n g_j\right) \sum_{j=k}^{n+1} \frac{(h_n - h_{j-1}) x_j + g_j}{w_{j+1,k+1}(h_k) w_{n,j-1}(h_n)}, \qquad 0 \le k \le n, \eqno(5.3)$$
and $L_{n,n+1}= 1$ for $n \ge 0$.
We can also write $L_{n,k}$ as 
$$L_{n,k}= \left( \prod_{j=k+1}^n g_j\right)\left(\sum_{j=k}^n \frac{x_j}{w_{n,j}(h_n)w_{j+1,k+1}(h_k)}+\sum_{j=\max(1,k)}^{n+1} \frac{g_j}{w_{n,j-1}(h_n) w_{j+1,k+1}(h_k)}\right). \eqno(5.4)$$

In particular
$$L_{n,n}= x_n + \frac{g_{n+1}}{h_n - h_{n+1}} -\frac{g_n}{h_{n-1} - h_n}, \eqno(5.5)$$ 
and 
$$ L_{n,n-1} = \frac{g_n}{h_{n-1} -h_n} \left(\frac{g_{n-1}}{h_{n-2} - h_n} - \frac{g_n}{h_{n-1} - h_n} + \frac{g_{n+1}}{h_{n-1}-h_{n+1}} +x_n - x_{n-1} \right). \eqno(5.6)$$

We will show in the next section that when the sequences $\{ h_k\}$, $\{x_k\}$, and $\{g_k\}$ are certain linearly recurrent sequences the matrix $L$ is tridiagonal and therefore the polynomial sequence $\{u_k(t): k \ge 0\}$ satisfies a three-term recurrence relation whose coefficients are the entries of $L$.

\section{The family of orthogonal polynomial sequences}
In this section we study the polynomial sequences obtained when the sequences $ h_k$ and $x_k$ satisfy a particular type of linear difference equation of third order and $g_k$ satisfies a linear difference equation of fifth order related with the equation satisfied by $ h_k$ and $x_k$. We will show that choosing appropriate initial values for the sequence $g_k$ the matrix $L$ becomes tridiagonal and thus the polynomial sequence $\{u_k(t): k \ge 0\}$ satisfies a three-term recurrence relation.

Let $z$ be a complex number and consider the difference equation
$$ s_{k+3} = z ( s_{k+2} -s_{k+1}) + s_k, \qquad k\ge -1. \eqno(6.1)$$
 We consider $s_0, s_1, s_2$ as the initial values of the sequence. In order to avoid undefined terms in some formulas it is convenient to define $s_{-1}=s_2 - z (s_1-s_0)$.
The characteristic polynomial of (6.1) is $ t^3 -z t^2 + z t -1$.
 The sum of its roots equals $z$, the product of the roots is equal to 1, and 1 is a root. 
 Therefore we can express the roots as 1, $q$, and $q^{-1}$ where $q$ is a nonzero complex number. If $z$ is real and $ -1 \le z \le 3$ then we can write $z= 1+ 2 \cos(\theta)$ and then it is easy to see that $q= \cos(\theta) + i \sin(\theta)$ and hence $q$ and $q^{-1}$ have modulus one.
If $z=3$ then 1 is a root of the characteristic polynomial with multiplicity three. If $z=-1$ then $-1$ is a double root. 

From now on we suppose that the sequences $h_k$ and $x_k$ are solutions of the difference equation (6.1) and that $h_k \ne h_n$ for $k \ne n$. It is easy to verify that the termwise (or Hadamard) product of two solutions of (6.1) satisfies the difference equation of order five
$$ s_{k+5}=(z^2-z-1) ( s_{k+4}- s_{k+1}) - (z-1)(z^2-z-1) ( s_{k+3} - s_{k+2}) +s_k, \qquad k\ge 0. \eqno(6.2)$$
The characteristic polynomial of this equation is $ (t^2-(z^2-2 z -1) t +1) ( t^3 -z t^2 + z t -1)$ and therefore its roots are $ 1,q,q^{-1}, q^2, q^{-2}$. 

Let $g_k$ be a sequence that satisfies (6.2) and also $g_0=0$ and $g_k\ne 0$ for $k \ge 1.$ The matrices $C$, $C^{-1}$, and $L$ are completely determined by $z$ and the initial values $ h_0,h_1,h_2$, $ x_0,x_1,x_2$, and $g_1,g_2,g_3,g_4$.  

Let us note that the entries of the matrices $C$, $C^{-1}$, and $L$ are functions of the differences $h_k- h_n$. We study next some properties of these differences.

Since the sequence $h_k$ satisfies (6.1) and $h_k- h_n \ne 0$ if $ k \ne n$, it is clear that 
$$\frac{z (h_{k+2} - h_{k+1})}{h_{k+3}- h_k} =1, \qquad k\ge 0, \eqno(6.3) $$
and then it is easy to see that
$$\frac{w_{3,0}(h_{k+1})}{w_{k+3,k}(h_1)}=-1, \qquad k\ge 2.$$ 
Shifting the sequence $h_k$ by adding $n$ to the indices in the previous equation we obtain 
$$\frac{w_{3+n,n}(h_{k+1+n})}{w_{k+3+n,k+n}(h_{1+n})}=-1,\qquad k \ge 2, \ n \ge 0,$$
and a change of variable yields a more symmetric expression
$$\frac{w_{3+n,n}(h_{k+1})}{w_{3+k,k}(h_{n+1})}=-1,\qquad n \ge k+2. \eqno(6.4)$$

For fixed $m >0$ the sequence $(\delta h)_{m,k} = h_k - h_{k+m} $ satisfies the difference equation
$$s_{k+2}= (z-1) s_{k+1} - s_k, \qquad k \ge 0, \eqno(6.5)$$
and has initial values 
 $(\delta h)_{m,0}= h_0- h_{m}$ and $(\delta h)_{m,1}= (h_1-h_{1+m})$. 
This is a simple three-term recurrence relation which is related with the recurrence satisfied by the Chebyshev families of orthogonal polynomials.

Let $p_k$, $r_k$, and $y_k$ be the solutions of (6.5) determined by the initial values $p_0=1$ and $p_1=z-1$; $r_0=1$ and $r_1=z$; $y_0=2$ and $y_1=z-1$. The sequence $p_k$ satisfies
$$p_k= U_k\left(\frac{z-1}{2}\right), \qquad k \ge 0, \eqno(6.6)$$
where the $ U_k(t)$ are the Chebyshev polynomials of the second kind. The sequence $r_k$ satisfies
$$r_k=W_k\left(\frac{z-1}{2}\right) =\frac{k!}{(1/2)_k} P_k^{1/2,-1/2} \left(\frac{z-1}{2}\right), $$
where the $W_k$ are the Chebyshev polynomials of the fourth kind, see \cite[(1.9)]{MasHan}, and the $P_k^{\alpha,\beta}(t)$ are Jacobi polynomials, and also satisfies  
 $r_{k+2}=z p_{k+1} - p_k$ for $k \ge 0$.
 The sequence $y_k$ satisfies 
$$y_k= T_k\left(\frac{z-1}{2}\right), \qquad k \ge 1, \eqno(6.7)$$
	where the $ T_k(t)$ are the Chebyshev polynomials of the first kind. The polynomials $y_k$ also satisfy $y_{k+2}=(z-1) p_{k+1} -2 p_k$ for $k \ge 0$.

If we put $p_{-1}=0$ and $p_{-2}= -1$ then the sequence $p_k$ and the shifted sequence $p_{k-1}$ are solutions of (6.5) for $k \ge -2 $. Since $ (\delta h)_{m,k}$ is also a solution of (6.5) it is easy to see that
$$ (\delta h)_{m,k+2}= (\delta h)_{m,1} p_{k+1} - (\delta h)_{m,0} p_k, \qquad k\ge 0, \ m\ge 1.\eqno(6.8)$$
Let us note that the termwise product of two solutions of (6.1) is a solution of (6.2) and every solution of (6.2) is also a solution of (6.1), with additional initial values. Therefore we can prove that 
$$ (\delta h)_{2 m,k}= (\delta h)_{2, k-1+m} p_{m-1}, \eqno(6.9)$$
and
$$ (\delta h)_{2 m+1,k}= (\delta h)_{1, k+m} r_m, \eqno(6.10)$$ 
 just verifying that the 4 initial values of each pair of sequences in the equations above coincide.
 
Therefore, the condition $h_k - h_n \ne 0$ for $k \ne n$ is satisfied if $z$ is not a root of any of the orthogonal polynomials $p_k$ and $r_k$ for $k\ge 0$ and $(\delta h)_{1, k}\ne 0$ and $(\delta h)_{2, k}\ne 0$ for $k \ge 0$. Note that $h_1\ne h_3$ is equivalent to $z \ne (h_0 -h_1)/(h_1-h_2)$.
 
\begin{theorem}
	Let $h_k$ and $x_k$ be sequences that satisfy the difference equation (6.1) and let $g_k$ be a sequence that satisfies (6.2) and has $g_0=0$. Then the matrix $L$, defined in (5.2), is tridiagonal if and only if 
$$g_3= z ((x_1-x_0)h_0 +(x_0-x_2) h_1 +(x_2-x_1) h_2 -g_1 + g_2), \eqno(6.11)$$
and 
$$g_4= (z-1) (z (x_2-x_0) h_0 + z(z+1) (x_1-x_2)h_1 -z ( z(x_1-x_2) -x_0+x_1)h_2 - (z+1)g_1 + z g_2) . \eqno(6.12) $$
\end{theorem}

The necessity of conditions (6.11) and (6.12) can be easily seen. Indeed, the conditions $L_{2,0}=0$ and $L_{3,1}=0$ are necessary for $L$ to be tridiagonal. Using the explicit expressions (5.3) or (5.4) for the entries of $L$ we see that equations (6.11) and (6.12) are equivalent to $L_{2,0}=0$ and $L_{3,1}=0$, respectively. The proof of the sufficiency of (6.11) and (6.12) is much harder. It will be given in Section 10.

Notice that the expressions for $g_3$ and $g_4$ in the theorem do not change if we interchange $x_j$ with $h_j$, for $ 0 \le j \le 2.$ Since $g_k$ satisfies the difference equation (6.2) we see that every $g_k$ is invariant under such interchange. Note also that, since $g_0=0$, when $g_3$ and $g_4$ satisfy (6.11) and (6.12), then $g_1$ and $g_2$ are free parameters that determine the sequence $g_k$.

From equation (5.3) we see that the entries of the matrix $L$ depend on the recurrent sequences $h_k, x_k, g_k$. Therefore it seems reasonable to expect that such entries satisfy some kind of recurrence relations. In what follows we present, without proofs, two such recurrence relations that can be used to prove Theorem 6.1, but since their proofs are quite long and messy, we give a different proof of Theorem 6.1 in Section 10. 

We present first a recurrence relation for three consecutive entries on a diagonal of $L$. The first step is to eliminate denominators of the entries of $L$. We also use properties (6.3) and (6.4) that allow us to simplify some quotients of polynomials $w(j,k,h_i)$.
Define
$$\tau(n,k) = \left(\prod_{j=k+1}^n \frac{1}{g_j} \right) w_{n+2,k+1}(h_k)\, w_{n,k+2}(h_n), \qquad 0 \le k \le n-2. \eqno(6.13)$$
Let us note that the denominator of $L_{n,k}$ is a divisor of $\tau(n,k)$. 

 The recurrence relation
$$ \tau(n+2,k+2) L_{n+2,k+2} - y_{n-k} \tau(n+1,k+1) L_{n+1,k+1} + \tau(n,k) L_{n,k}=0, \qquad 0 \le k \le n-2, \eqno(6.14)$$  
holds, where $y_j$ is the sequence of orthogonal polynomials in $z$ related with the monic Chebyshev polynomials of the first kind that satisfies (6.7).

 We present next a recurrence relation on $L$ similar to the recurrences satisfied by the entries of the matrix of Stirling numbers of the second kind and other matrices of generalized binomial coefficients.

 Define the functions
 $$\epsilon_1(n,k)= w_{n+2,k+1}(h_k)\, w_{n,k+2}(h_n), $$
 $$\epsilon_2(n,k)=g_{n+1}\, w_{k+1,k}(h_{k+1})\, w_{n+3,k+1}(h_k) \frac{w_{n+1,k+2}(h_{n+1})}{w_{n+1,n-1}(h_{n+1})}, $$
 $$\epsilon_3(n,k)= g_k \, w_{n,n-1}(h_n)\, w_{n+2,k+2}(h_{k-1})\, w_{n,k+1}(h_n). \eqno(6.15)$$

The entries of $L$ satisfy the recurrence relation
$$p_{n-k-2} \epsilon_1(n,k) L_{n,k} -\epsilon_2(n-1,k) L_{n-1,k} + \epsilon_3(n,k+1) L_{n,k+1}=0, \qquad 0 \le k \le n-2, \eqno(6.16)$$
where $p_j$ is the sequence of orthogonal polynomials of $z$ that satisfies (6.6).

From the recurrence relations (6.14) and (6.16) we can see that if $L_{2,0}=0$ and $L_{3,1}=0$ then $L_{m+k,k}=0$ for $m\ge 2$ and $k\ge 0$. That is, $L$ is tridiagonal. 

\section{The class $\cH_q$}
In this case the roots of the polynomial $ t^3 -z t^2 + z t -1$ are $1, q, q^{-1}$, with $q \ne 1$ and $q \ne -1$. Therefore the roots are distinct and then the general solution of the difference equation (6.1) can be expressed as
$$s_k = \lambda_0 + \lambda_1 q^k +\lambda_2 q^{-k}, \qquad k \ge 0, \eqno(7.1)$$
where the coefficients $\lambda_j$ are complex numbers. Since the sequences $h_k$ and $x_k$ are solutions of (6.1) they can be expressed as 
$$h_k= a_0 + a_1 q^k + a_2 q^{-k}, \qquad k \ge 0, \eqno(7.2)$$
 where
$$ a_0= \frac{(q^2+1) h_1 -q (h_0+h_2)}{(q-1)^2},$$ 
$$ a_1=- \frac{q(h_1 -h_2) +h_1-h_0}{(q+1) (q-1)^2},$$ 
$$ a_2=\frac{q^2 (q (h_0-h_1) + h_2 -h_1)}{(q+1) (q-1)^2}, \eqno(7.3)$$
and 
$$x_k= b_0 + b_1 q^k + b_2 q^{-k}, \qquad k \ge 0, \eqno(7.4)$$
where 
$$ b_0= \frac{(q^2+1) x_1 -q (x_0+x_2)}{(q-1)^2},$$ 
$$ b_1=- \frac{q(x_1 -x_2) +x_1-x_0}{(q+1) (q-1)^2},$$ 
$$ b_2=\frac{q^2 (q (x_0-x_1) + x_2 -x_1)}{(q+1) (q-1)^2}. \eqno(7.5)$$
Let us note that a change in $b_0$ corresponds to a translation of the Newton polynomials $v_k(t)$ and therefore to a translation of the polynomials $u_k(t)$, and a change in $a_0$ corresponds to a translation of the $h_k$.

The characteristic roots of the difference equation (6.2) are $1,q,q^{-1}, q^2, q^{-2}$ and therefore the general solution of (6.2) has the form
 $$ s_k= \lambda_0 + \lambda_1 q^k +\lambda_2 q^{-k} + \lambda_3 q^{2k} + \lambda_4 q^{-2k}, \qquad k \ge 0, \eqno(7.6)$$
where the coefficients $\lambda_j$ are complex numbers. Since the sequence $g_k$ is a solution of (6.2) it can be expressed as 
$$g_k= d_0 + d_1 q^k + d_2 q^{-k} + d_3 q^{2k} + d_4 q^{-2k} , \qquad k \ge 0. \eqno(7.7)$$
Since $g_0=0$ and $g_3$ and $g_4$ are given by (6.11) and (6.12) respectively, we obtain
$$d_0= -( a_2 b_2 q + d_1 + d_2 + a_1 b_1 q^{-1}),\qquad d_3=a_1 b_1 q^{-1}, \qquad d_4= a_2 b_2 q, \qquad \eqno(7.8)$$
and $d_1$ and $d_2$ are arbitrary parameters. Therefore we have
$$g_k= (q^k -1) d_1 + (q^{-k} -1) d_2 + q^{-1} (q^{2k} -1) a_1 b_1 + q (q^{-2 k}-1) a_2 b_2, \qquad k \ge 0.\eqno(7.9)$$

The coefficients of the three-term recurrence relation are given in equations (5.5) and (5.6). We can express them in terms of the parameters $a_0,a_1,a_2,b_0,b_1,b_2, d_1,d_2$ by using equations (7.2), (7.4), and (7.7). The entries $\alpha_k = L_{k,k-1}$, for $k \ge 1$, are given by
$$ \alpha_k = \frac{ (q^{k-1} a_1 + a_2) (q^{2k} a_1 b_2 - a_2 b_1) + q^k (q^{k-1} a_1 d_2 - a_2 d_1) }{(q^{2 k} a_1-a_2) (q^{2 k -1} a_1-a_2)^2 (q^{2 k-2} a_1-a_2) } \times \qquad \qquad \qquad \qquad \hfill $$ 
$$ \qquad \qquad \qquad (q^k-1) (q^{k-1} a_1 - a_2)
 \ ( (q^k+1) (q^{2k-2} a_1 b_1 - a_2 b_2)+ q^{k-1} ( q^k d_1 -d_2)). \eqno(7.10)$$
Let $\beta_k=L_{k,k}$ and let $\sigma_k=\beta_0+\beta_1+ \cdots +\beta_k$, for $k \ge 0$. Then we obtain 
$$\sigma_k= (k+1) b_0 +\left( \frac{q^{k+1}-1}{q-1}\right) \frac{ (q^k a_1 + a_2) ( q b_2 -b_1) - q^{k+1} d_1 +d_2}{q^{2k+1} a_1 - a_2}, \eqno(7.11)$$

The polynomials $u_k(t)$ are invariant under the maps  
$$(q,q^{-1},a_1,a_2,b_1,b_2, d_1,d_2) \rightarrow (q^{-1}, q, a_2,a_1,b_2,b_1,d_2,d_1) $$ 
and 
$$ (a_1,a_2,d_1,d_2) \rightarrow r (a_1,a_2,d_1,d_2), \qquad r \ne 0.$$
This fact can be seen using equations (4.2), (4.3), (7.2), (7.7), and (7.8). 
This clearly implies that the recurrence coefficients $\alpha_k$ and $\beta_k$ are also invariant under such transformations.

In some concrete cases these properties, combined with symmetries of the parameters that determine the concrete case, are the reason why in some cases we can find different sets of values of $ a_1,a_2, b_0,b_1,b_2, d_1,d_2$ that produce the same family of polynomial sequences. Let us note that multiplication of $b_0, b_1,b_2, d_1, d_2$ by a nonzero constant produces a rescaling of the $x_k$ which is equivalent to a rescaling of the variable $t$ in the Newton type basis $\{v_k(t)\}$ and also in the sequence $\{u_k(t)\}$.

From (7.10) we see that in order to have $\alpha_k \ne 0$, for $k \ge 1$, the parameter $q$ must not be a root of 1. Let us note that if $q^k a_1 -a_2 \ne 0$ for $ k\ge 1$ then $\alpha_k$ and $\beta_k$ are well defined. Note also that $\alpha_k$ is independent of $a_0$ and $b_0$, and $\beta_k$ is independent of $a_0$.  

Equation (7.2) gives us 
$$ h_k - h_j = (q^{k} -q^{j}) a_1 +(q^{-k} - q^{-j}) a_2 = (q^{-j} - q^{-k}) (q^{k+j} a_1 - a_2), \qquad k \ne j, \eqno(7.12)$$
and therefore we must have $ q^n a_1 \ne a_2$ for $ n \ge 1$.
With equations (7.12) and (7.9) we can express the entries of the matrices $C$ and $C^{-1}$ in terms of the parameters $a_0,a_1,a_2,b_0,b_1,b_2, d_1,d_2$. For example, the quotients
$$\frac{g_k}{h_0 - h_k}=\frac{ d_2 - q^k d_1 -q^{k-1} (1+ q^k ) a_1 b_1 +q (1+q^{-k}) a_2 b_2}{q^k a_1 - a_2}, \qquad k\ge 1, \eqno(7.13)$$  
 are the factors of the generalized moments $m_n=\hat c_{n,0}$, given by equation (4.9). 

In the generalized difference equation (4.1) the operators $\gamma$ and $\phi$ are defined with respect to the basis $\{v_k(t): k \ge 0 \}$. Using the matrices $V$ and $V^{-1}$ we can transform equation (4.1) and obtain an equation with respect to the basis $\{t^k: k\ge 0\}$. For example, if $h_0=0$, $b_0=0$, and $b_2=0$ we obtain the $q$-difference equation of second order
$$ \left( \frac{(q-1)^2}{q^2} f_2 (t) \cD_{1/q} \cD_q + \frac{q-1}{q} f_1(t) \cD_q \right) u_k(t) = h_k u_k(t), \qquad k \ge 0, \eqno(7.14)$$
where
$$f_1(t)= (q a_1-a_2) t +q d_1 + a_1 b_1 + a_2 b_1 - d_2, \qquad f_2(t)= a_2 t^2 +(d_2-a_2 b_1) t - b_1 d_2, $$
$h_k$ is given by $(7.2)$ and $ \cD_q$ and $ \cD_{1/q}$ are the usual $q$-difference and $q^{-1}$-difference operators. There is a similar equation in the case with $h_0=0$, $b_0=0$, and $b_1=0$. If $b_1$ and $b_2$ are both nonzero, a second order $q$-difference equation as for the Askey-Wilson polynomials is obtained.

If we know all the ingredients of the operator in the left-hand side of (7.14) we can verify that (7.14) is equivalent to (4.1) by applying the operator to $v_k(t)$. 
In the general case the matrix representation of the generalized difference operator with respect to the basis of monomials is not a banded matrix.

The class $\cH_q$ contains all the families of basic hypergeometric orthogonal polynomial sequences in the $q$-Askey scheme \cite{Hyp}. The coefficients of the normalized three-term recurrence relation for each of the families listed in Chapter 14 of \cite{Hyp} are obtained by direct substitution of appropriate values of the parameters $ a_1,a_2,b_0, b_1, b_2, d_1,d_2$, without taking limits. The only case that requires taking a limit is the family of continuous $q$-Hermite polynomials. Some families correspond to several sets of values of the parameters. We give next some examples.

The Askey-Wilson polynomials are obtained with
$$ a_1= a b c d q^{-1} a_2, \qquad b_0=0, \qquad b_1=a/2, \qquad b_2=a^{-1}/2, $$
$$ d_1=- a ( abcd +q (b c + b d + c d)) q^{-2} a_2/2, \qquad d_2=- ( (b + c + d) + q a^{-1} ) a_2/2, \eqno(7.15)$$
where $a_2$ is an arbitrary nonzero number and $a,b,c,d$ are the parameters used in \cite[eq. 14.1.5]{Hyp}. We obtain other representations for this family interchanging $a$ with $b$, $c$, or $d$ in the expressions for the parameters $a_j, b_j, d_j$ in (7.15).

The $q$-Racah polynomials are obtained with
$$ a_1=\alpha \beta q a_2, \qquad b_0=0, \qquad b_1=\gamma \delta q, \qquad b_2=1, $$
$$ d_1= -q ( \alpha \beta \gamma \delta + \alpha \beta \delta + \beta \gamma \delta +
\alpha \gamma ) a_2 , \quad d_2=-q ( \beta \delta + \alpha + \gamma +1) a_2, \eqno(7.16)$$
where $a_2$ is an arbitrary nonzero number and $ \alpha, \beta, \gamma, \delta $ are the parameters 
used in \cite[eq. 14.2.4]{Hyp} and must satisfy 
$$ \alpha q= q^{-N}, \qquad \hbox{or}\ \ \beta \delta q=q^{-N}, \qquad \hbox{or}\  \gamma q = q^{-N}, $$
 where $N$ is a positive integer.

The continuous dual $q$-Hahn polynomials are obtained with
$$ a_1=0,\qquad b_0=0, \qquad b_1=\frac{a}{2},\qquad  b_2=\frac{1}{2 a}, $$
$$ d_1= \frac{-a b c a_2}{2 q}, \qquad d_2=\frac{-(a b + ac + q) a_2}{2 a}, \eqno(7.17)$$
where $a_2$ is an arbitrary nonzero number and $a,b,c$ are the parameters used in \cite[eq. 14.3.5]{Hyp}.

The Al-Salam-Chihara polynomials are obtained with
$$ a_1=0,\qquad b_0=0, \qquad b_1=a,\qquad b_2=a^{-1}, $$
$$ d_1= 0, \qquad d_2= - a_2 ( b + q a^{-1}), \eqno(7.18)$$
where $a_2$ is any nonzero number and $ a, b$ are the parameters in \cite[eq. 14.8.5]{Hyp}.

For the big $q$-Jacobi polynomials we have
$$ a_1=a b q a_2, \qquad b_0=0,\qquad b_1=0, \qquad b_2=1, $$
$$ d_1= -a c q a_2, \qquad d_2= -(1+ a +c ) q a_2, \eqno(7.19) $$
where $a_2$ is any nonzero number and $a,b,c$ are the parameters in \cite[eq. 14.5.4]{Hyp}.
Another set of values of the parameters that gives the same family of big $q$-Jacobi polynomials is 
$$ a_1=a b q a_2, \qquad b_0=0,\qquad b_1=a q,\qquad b_2=0, $$
$$ d_1= - (a b + b + c) a q a_2 , \qquad d_2=-c q a_2. \eqno(7.20)$$
Note that the Newton bases for the two sets of values are different.

For the $q$-Meixner polynomials we have
$$ b_0=0, \qquad b_1=b q , \qquad b_2=0, \qquad a_2=0,$$
$$ d_1=a_1 ( b c - b -1), \qquad d_2= c a_1, \eqno(7.21)$$
 where $a_1$ is an arbitrary nonzero number and $b,c$ are the parameters in \cite[eq. 14.13.4]{Hyp}.
 This family is also obtained with $(b_0=0, b_1=0, b_2=1, a_2=0)$ and with $( (b_0=0, b_1=0, b_2=- c b, a_2=0)$, taking suitable values for $d_1$ and $d_2$ for each set of values of the parameters. Note that the corresponding sequences of nodes $x_k$ are different in the three cases.

\section{The class $\cH_1$}
In this section we consider the family of orthogonal polynomial sequences obtained when we take $z=3$. In this case 
 the characteristic polynomial of the difference equation (6.1) is $ t^3 - 3 t^2 + 3 t -1 = (t-1)^3$ and therefore the general solution of (6.1) has the form
 $$ s_k = \lambda_0 + \lambda_1 k + \lambda_2 k ( k-1) , \qquad k\ge 0, \eqno(8.1)$$
 where the coefficients $\lambda_j$ are arbitrary complex numbers. 

Since the sequences $h_k$ and $x_k$ are solutions of (6.1) we can write them as follows.
$$h_k= a_0 + a_1 k + a_2 k ( k-1), \qquad k \ge 0, \eqno(8.2)$$
where $ a_0=h_0$, $a_1=h_1 - h_0$, and $ a_2= (h_0 -2 h_1 + h_2)/2$, and
$$ x_k= b_0 + b_1 k + b_2 k ( k-1), \qquad k \ge 0, \eqno(8.3)$$
where $ b_0=x_0$,  $b_1=x_1 - x_0$, and $ b_2= (x_0 -2 x_1 + x_2)/2$.

The characteristic polynomial of the difference equation (6.2) is in this case $(t-1)^5$ and therefore the
sequence $g_k$, which is a solution of (6.2), with $g_0=0$ and $g_3$ and $g_4$ given by (6.11) and (6.12) respectively, can be expressed as
$$ g_k= d_1 k + d_2 k (k-1) + d_3 k (k-1) (k-2) + d_4 k (k-1) (k-2)(k-3), \qquad k \ge 0, \eqno(8.4)$$
where 
$$d_1= g_1, \quad d_2= \frac{g_2}{2} -g_1, \quad d_3=a_1 b_2 + a_2 b_1 + 2 a_2 b_2, \quad d_4=a_2 b_2, \eqno(8.5)$$
and $g_1$ and $g_2$ are any nonzero numbers.

From (8.2) we obtain 
$$ h_m - h_n = (m-n) ( a_1+ (m -1 +n) a_2), \qquad m \ge 0,\ n \ge 0, $$
and therefore the condition $h_m \ne h_n$, for $m \ne n$, is satisfied if $a_1 + k a_2 \ne 0$ for $ k \ge 0$.

The coefficients of the three-term recurrence relation are in this case given by
$$\alpha_k= \frac{g_k \ ( a_1 +(k-2) a_2) }{(a_1+(2 k-1) a_2) (a_1 +(2k-2) a_2)^2 (a_1+ (2k-3) a_2)} \times \qquad \qquad \qquad \qquad $$ 
$$ \qquad \qquad (( a_1+(k-1) a_2) (a_2 b_2 k (k-1) +(a_1 b_2 -a_2 b_1) k -a_1 b_1 -a_1 b_2 + d_2 ) - d_1 a_2 ), \eqno(8.6)$$
 where $g_k$ is given by (8.4) and 
$$ \sigma_k= (k+1) \left( b_0 - \frac{d_1 + e_1 k + e_2 k (k-1) + e_3 k (k-1) (k-2) }{a_1 + 2 k a_2 }\right) ,\qquad k \ge 0, \eqno(8.7)$$
where $e_1= -(1/2) (a_1 b_1 + 2 a_2 b_1 - 2 d_2),$ $ e_2= (2/3) b_2 ( a_1 + a_2),$ $ e_3=(1/3) a_2 b_2 $.  
Recall that $\beta_0=\sigma_0$ and $\beta_k= \sigma_k - \sigma_{k-1}$ for $k \ge 1$.

The recurrence coefficients $\alpha_k$ and $\beta_k$ are invariant under the map
$$(a_1,a_2,d_1,d_2)\rightarrow r (a_1,a_2,d_1,d_2),\qquad r \ne 0.$$

The quotients 
$$ \frac{g_k}{ h_0 - h_k} =- \frac{d_1 + d_2 (k-1)+2 (a_1 b_2 + a_2 b_1 +2 a_2 b_2 )\binom{k-1}{2} +6 a_2 b_2 \binom{k-1}{3}}{a_1 + (k-1) a_2}, \qquad k\ge 1, \eqno(8.8)$$
are the factors of the generalized moments $m_n=\hat c_{n,0}$ given by (4.9).
The entries of the matrices $C$ and $C^{-1}$ are products of terms similar to (8.8).

If $b_0$, $b_1$, and $b_2$ are equal to zero then the basis $\{v_k(t)\}$ is the standard basis of monomials $\{t^k\}$ and in such case the generalized difference equation (4.1) can be written in the form
$$ ((d_2 t + a_2 t^2) D^2 +( d_1 + a_1 t) D + a_0) u_k(t) = h_k u_k(t), \qquad k \ge 0, \eqno(8.9)$$
where $D$ denotes differentiation with respect to $t$ and $h_k$ is given by (8.2).

If we take $b_1=0$ and $b_2=0$ then the basis $\{v_k(t)\}$ becomes the basis of translated powers $\{(t-b_0)^k: k \ge 0\}$ and the matrices $V$ and $V^{-1}$ are generalized Pascal matrices. Transforming the generalized difference equation (4.1) to the basis of monomials we obtain the equation 
$$ ( (b_0 (a_2 b_0 -d_2) + ( d_2 - 2 a_2 b_0) t + a_2 t^2) D^2 + (d_1 -a_1 b_0 + a_1 t) D + a_0) u_k(t)=h_k u_k(t), \eqno(8.10)$$
where $D$ denotes differentiation with respect to $t$. If we put $b_0=0$ in (8.10) we obtain (8.9).

If $b_2=0$ and $b_1\ne 0$ then the basis $\{v_k(t)\}$ is the Newton basis associated with the sequence $b_0 + b_1 k$, for $k \ge 0$. 
In this case the generalized difference operator $\cD$ of (4.1) can be expressed as 
$$ \cD= a_0 I +((a_1-2 a_2) t + \epsilon_1) \Delta_{b_1} + (a_2 t + \epsilon_2 ) \Delta_{b_1}^2 \, E , \eqno(8.11)$$
where $\Delta_{b_1}$ is the difference operator defined by 
$$ \Delta_{b_1} v_k(t) =\frac{v_k(t+b_1) - v_k(t)}{b_1}= k v_{k-1}(t), \qquad k \ge 0, \eqno(8.12)$$
Here $E$ is the forward shift operator defined by
 $$ E v_k(t) = v_{k+1}(t), \qquad k \ge 0, \eqno(8.13)$$
 and the coefficients $\epsilon_1, \epsilon_2 $ are given by
$$ \epsilon_1=d_1 - 2 d_2 -a_1 b_0 +2 a_1 b_1 + 2 a_2 b_0 + 2 a_2 b_1, \qquad \epsilon_2=d_2 = a_1 b_1 -a_2 b_0 - a_2 b_1. $$ 

 The class $\cH_1$ contains all the families of hypergeometric orthogonal polynomial sequences in the Askey scheme. The coefficients of the normalized three-term recurrence relation for each family listed in Chapter 9 of \cite{Hyp} are obtained by giving appropriate values to our parameters $ a_1,a_2, b_0,b_1,b_2,d_1,d_2$, without taking limits. The only case that requires limits is the family of Hermite polynomials. We give next some examples.

The Wilson polynomials are obtained with
$$ a_1=(a+b+c+d) a_2, \qquad b_0=-a^2, \qquad b_1=-2 a -1, \qquad b_2= -1, $$
$$ d_1= -a_2 (a+b) (a+c) (a+d), \quad d_2=-a_2 ((2 a +1) ( a + b + c + d +1) + a^2 + bc+ bd+ cd),\eqno(8.14)$$
where $a_2$ is any nonzero number and $a,b,c,d$ are the parameters in \cite[eq. 9.1.5]{Hyp}.

The Racah polynomials, in the case $\alpha +1=-N$, are obtained when
$$ a_1=(1 + \beta -N) a_2, \qquad b_0=0,\qquad b_1=2 +\gamma + \delta, \qquad b_2=1,$$
$$ d_1=- N ( 1 +\gamma) ( 1 + \beta + \delta) a_2, \qquad d_2= -a_2 ((3+ \beta + \gamma+\delta) N - ( 2 +\gamma)( 2 +\delta + \gamma)), \eqno(8.15)$$
where $a_2$ is any nonzero number and $ \alpha \beta, \gamma, \delta, N$ are the parameters in \cite[eq. 9.2.4]{Hyp}.

For the continuous dual Hahn polynomials we have
$$ a_2=0,\qquad b_0=-a^2,\qquad b_1=-2 a -1,\qquad b_2=-1,$$
$$d_1=-a_1 (a+c) (a+b), \qquad d_2=-a_1 ( 2 a + b + c +1), \eqno(8.16)$$
where $a_1$ is any nonzero number and $ a,b,c$ are the parameters in \cite[eq. 9.3.5]{Hyp}.

For the continuous Hahn polynomials we have
$$ a_1=a_2 (a+b+c+d), \qquad b_0=i b, \qquad b_1=i,\qquad b_2=0,$$
$$ d_1= i a_2 ( b^2+ b c + b d + c d), \qquad d_2= i a_2 ( 2 b + c + d + 1), \eqno(8.17)$$
where $a_2$ is any nonzero number and $ a,b,c,d$ are the parameters in \cite[eq. 9.4.4]{Hyp}.

 The Meixner-Pollaczek polynomials are obtained with
 $$ a_2=0, \qquad b_0= -i \lambda, \qquad b_1= -i,\qquad b_2= 0, $$
 $$ d_1=a_1 \lambda \left( \frac{\cos(\phi)}{\sin(\phi)} -i\right), \qquad d_2= \frac{a_1}{2} \left( \frac{\cos(\phi)}{\sin(\phi)} -i\right), \eqno(8.18)$$
where $a_1$ is any nonzero number and $ \lambda, \phi $ are the parameters in \cite[eq. 9.7.4]{Hyp}.

For the Jacobi polynomials we have
$$ a_1= a_2 ( 2 + \alpha + \beta), \qquad b_0=1,\qquad b_1=0, \qquad b_2=0, $$
$$ d_1= 2 a_2 ( \alpha +1), \qquad d_2=2 a_2, \eqno(8.19)$$
here $a_2$ is any nonzero number and $\alpha, \beta$ are the parameters in \cite[eq.9.8.5]{Hyp}.

The Bessel polynomials are obtained with
$$ a_1=( a+2) a_2, \qquad b_0=0, \qquad b_1=0,\qquad b_2=0,\qquad 
d_1= 2 a_2, \qquad d_2=0, \eqno(8.20)$$
where $a_2$ is any nonzero number and $a$ is the parameter in \cite[eq. 9.13.4]{Hyp}.

\section{The class $\cH_{-1}$}
We consider now the class of orthogonal polynomial sequences obtained when we take $z=-1$. In this case the characteristic polynomial of the difference equation (6.1) is $t^3 + t^2 - t -1= ( t-1) (t+1)^2$ and its roots are $ 1,-1,-1$. Therefore the general solution of (6.1) has the form
$$ s_k = \lambda_0 + \lambda_1 (-1)^k + \lambda_2 k (-1)^k, \qquad k\ge 0, \eqno(9.1)$$
 where the $\lambda_j$ are arbitrary numbers.
 
 We can write the sequences $h_k$ and $x_k$, which are solutions of (6.1), as 
 $$ h_k = a_0 + a_1 (-1)^k + 2 a_2 k (-1)^k, \qquad k\ge 0, \eqno(9.2)$$
 and 
 $$ x_k= b_0 + b_1 (-1)^k + 2 b_2 k (-1)^k, \qquad k \ge 0. \eqno(9.3) $$
 We include the factor $2$ in the terms with $a_2$ and $b_2$ in order to simplify the notation.

When $z=-1$ the difference equation (6.2) has characteristic roots $ 1,1,1,-1,-1$ and then we can write
$$ g_k = d_0 + d_1 (-1)^k + 2 d_2 k (-1)^k + 2 d_3 k + 2 d_4 k (k-1), \qquad k \ge 0. \eqno(9.4)$$
Since $g_0=0$ and $g_3$ and $g_4$ must satisfy (6.11) and (6.12) respectively, we obtain
$$ d_0= - d_1, \qquad d_3= -a_1 b_2 - a_2 b_1, \qquad d_4= - 2 a_2 b_2. \eqno(9.5)$$
The coefficients $d_1$ and $d_2$ are arbitrary numbers.

Equation (9.2) gives us 
$$ h_n - h_m = ((-1)^n - (-1)^m ) a_1 + 2( (-1)^n n - (-1)^m m ) a_2, $$ 
and this shows that the condition $h_n - h_m \ne 0$ for $ n \ne m$ is satisfied if $a_2$ is not zero and $a_1 - n a_2 \ne 0$ for $n \ge 1$.

The initial terms of the sequence $h_k$ are 
$$ a_0 + a_1, a_0 -a_1 - 2 a_2, a_0 +a_1+4 a_2, a_0-a_1 - 6 a_2, \ldots .$$
Note that the terms with even indices and the terms with odd indices behave differently. The sequences $x_k$ and $g_k$ have a similar property. 

The coefficients of the three-term recurrence relation also have different expressions for even and odd indices.
For $n$ even we have
$$\alpha_n=- \frac{n ( a_1 + (n-1) a_2) ( a_1 b_2 + a_2 b_1 - d_2 + 2 (n -1) a_2 b_2) ( 2n a_2 b_2 + a_1 b_2 - a_2 b_1 + d_2)}{ a_2 ( a_1 + (2n -1) a_2)^2} . \eqno(9.6)$$
For $ n$ odd we have
$$\alpha_n= \frac{-d_1 - n (a_1 b_2 + a_2 b_1 + d_2) - 2 n ( n-1) a_2 b_2}{a_2 (a_1 + (2 n -1) a_2)^2} \times \qquad \qquad \qquad \qquad $$
$$ ((a_1 (a_1 b_2 - a_2 b_1 - d_2) + a_2 d_1 + a_2 ( (3 n -1) a_1 b_2 -(n-1)( a_2 b_1 +d_2) + 2 n ( n-1) a_2 b_2))). \eqno(9.7)$$
Note that, in both cases, $\alpha_n$ is a rational function of $n$ with numerator of degree 4 and denominator of degree 2.  

For the coefficients $\sigma_n = \beta_0+ \beta_1+\cdots+ \beta_n$ we have, if $n$ is even 
$$ \sigma_n= (n+1) b_0 + \frac{n (a_2 b_1 - a_2 b_2 -d_2) +a_1 b_1 - a_1 b_2 -d_1 -d_2}{a_1+( 2 n +1) a_2},\eqno(9.8)$$
and if $n$ is odd then
$$ \sigma_n= ( n+1) \left( b_0+ \frac{a_2 b_1 -a_2 b_2 - d_2)}{ a_1 + (2 n +1)a_2} \right). \eqno(9.9)$$

 The generalized moments, given by (4.9), are products of the quotients
 $$ \frac{g_k}{h_0 - h_k} = \frac{ 4 k^2 a_2 b_2 - (-1)^{k} (d_1 +2 k d_2) - 2 k ( 2 a_2 b_2 - a_1 b_2 - a_2 b_1)+ d_1}{ (-1)^k (a_1 + 2 k a_2) - a_1}, \qquad k \ge 1. \eqno(9.10)$$

 We consider next the generalized difference equation (4.1) in a particular case.
 If we take $b_0=0$ and $b_2=0$ then the matrix representation $B$ of the difference operator in (4.1) with respect to the basis of monomials is a sum of Kronecker products of infinite matrices and $2 \times 2$ matrices. That is 
 $$B= I \otimes B_0 + D \otimes B_1 + (D S\tr) \otimes B_2, \eqno(9.11)$$ 
 where $\otimes $ denotes the Kronecker product and 
 $$ B_0=\left[ \begin{matrix} a_0 + a_1 & 0 \cr 2 a_1 b_1 - 2 d_1 - 2 d_2 & a_0 - a_1 - 2 a_2 \cr \end{matrix} \right], \qquad B_1= \left[ \begin{matrix} -4 b_1 d_2 & -4 a_2 b_1 + 4 d_2 \cr 0 & 4 b_1 d_2 \cr \end {matrix} \right], $$
	 and 
$$ B_2=\left[ \begin{matrix} 4 a_2 & 0 \cr 4 a_2 b_1 - 4 d_2 - 4 a_2 & 0 \cr \end{matrix} \right]. $$

	Let $A= C V$ be the matrix of coefficients of the sequence $\{u_k(t)\}$ with respect to the monomial basis. Then the difference equation (4.1) is in this case equivalent to the equation $ A B = H A$. The operators that appear here seem to be related with Dunkl operators. We have not studied particular families in $\cH_{-1}$.

Some orthogonal polynomial sequences in the class $\cH_{-1}$ have been studied in the recent papers \cite{Vinet1}, \cite{Vinet2}, and \cite{Vinet3}. Such sequences are obtained by taking limits as $q$ goes to $-1$ and are connected with Hahn and big and little $q$-Jacobi polynomials.  
	
\section{ The proof of Theorem 6.1}

The proof of Theorem 6.1 will be obtained in three steps. In the first step we show that $ g_k= x_{k-1}(h_k - h_0) + d_k$ satisfies (6.2), (6.11), (6.12), and $g_0=0$ if and only if $d_k$ satisfies (6.1) and $d_0=0$. In the second step we show that a certain matrix $R$, constructed with the sequences $h_k$ and $d_k$, is tridiagonal if and only if $d_k$ satisfies (6.1) and $d_0=0$. In the last step we show that $L$ is tridiagonal if and only if $R$ is tridiagonal. Combining these results we obtain a proof of Theorem 6.1.

 In this Section we assume that  $h_k$ and $x_k$ are sequences that satisfy the difference equation (6.1) for $k \ge -1$. (Note that sequences given for $k\ge 0$ and satisfying (6.1) are immediately extended to such sequences for $k \ge -1$.)

   Let $g_k$ and $d_k$ be sequences for $k \ge 0$ which are related by $g_k= x_{k-1}(h_k - h_0) + d_k$. It is easy to verify that the termwise product of two sequences that satisfy the difference equation (6.1) is a sequence that satisfies (6.2). Hence $g_k$ satisfies (6.2) and $g_0=0$ if and only if $d_k$ satisfies (6.2) and $d_0=0$.
Since $ d_1= g_1 -x_0 ( h_1-h_0) $ and $d_2= g_2- x_1 ( h_2 - h_0) $, substitution in $g_3$ and $g_4$ gives us 
$$g_3= z ((x_1-x_0)h_0 +(x_0-x_2) h_1 +(x_2-x_1) h_2 -g_1 + g_2). \eqno(10.1)$$
and 
$$g_4= (z-1) (z (x_2-x_0) h_0 + z(z+1) (x_1-x_2)h_1 -z ( z(x_1-x_2) -x_0+x_1)h_2 - (z+1)g_1 + z g_2) . \eqno(10.2) $$
These equations are the same as (6.11) and (6.12), respectively, which are equivalent to $L_{2,0}=0$ and $L_{3,1}=0$.

From (10.1) and (10.2) we get 
$$g_3-z((x_1-x_0)h_0 + (x_0-x_2)h_1 +(x_2-x_1) h_2 + g_2 -g_1)= d_3-z(d_2-d_1)$$
and
$$g_4-(z-1)(z(x_2-x_0)h_0+z(z+1)(x_1-x_2)h_1-z(z(x_1-x_2)-x_0+x_1)h_2-(z+1)g_1 +z g_2)$$
$$= (d_4-z(d_3-d_2)-d_1)+z(d_3-z(d_2-d_1)).$$
(Note that the characteristic polynomial of (6.1) divides the one of (6.2), and thus any solution of (6.1) is also a solution of (6.2).)
Hence $g_k$ satisfies (6.2), (10.1), and (10.2) if and only if $d_k$ satisfies (6.1). This completes the first step.

For the second step, let $K$ be the infinite lower triangular matrix obtained by putting $g_k=1$, for all $k \ge 1$, in the matrix $C$, defined in Section 4. Note that $K$ depends only on the sequence $h_k$. 
 From equations (4.3), (4.6), and (4.7) we see that the $(n,j)$ entry of $K$ is $1/w_{n,j}(h_n)$ and the $(j,k)$ entry of $K^{-1}$ is $1/w_{j+1,k+1}(h_k)$. 

 Let $\{d_k:k \ge 0\}$ be a sequence that satisfies (6.1) and has $d_0=0$, and let $\Delta$ denote the infinite diagonal matrix whose $(k,k)$ entry is $d_{k+1}$. Recall that $S$ is the shift operador defined in (3.5). We define the matrix $R= K \, \Delta \, S\tr \, K^{-1}$.
 
 \begin{theorem}
	 The matrix $R$ is tridiagonal if and only if the sequence $d_k$ satisfies (6.1).
 \end{theorem}
 
{\it Proof:} The expressions for the entries of $K$ and $K^{-1}$ give us  
	 $$R_{n,k}= \sum_{j=k}^{n+1} \frac{d_j}{w_{n,j-1}(h_n)\, w_{j+1,k+1}(h_k) }. \eqno(10.3)$$
	Since $ h_{-1}=h_2-z(h_1-h_0)$ there are no undefined terms in the denominator when $j=0$. 
	We must show that $R_{n,k}=0$ for $n \ge 2$ and $0 \le k \le n-2$.

	It is convenient to define 
	 $$\hat R_{n,k}= w_{n+2,k+1}(h_k) \, R_{n,k}, \qquad n \ge 2,\ 0 \le k \le n-2.$$
Since the $h_k$ are pairwise distinct $w_{n+2,k+1}(h_k) \ne 0$, and then $ R_{n,k}=0$ if and only if $\hat R_{n,k}= 0$.
A simplification gives us 
	 $$\hat R_{n,k}= \sum_{j=k}^{n+1} \frac{w_{n+2, j+1}(h_k)\, d_j}{w_{n,j-1}(h_n)}. $$
Define the polynomials 
	 $$y_{n,k}(t) = \sum_{j=k}^{n+1} \frac{w_{n+2, j+1}(h_k)\, t^{j-k}}{w_{n,j-1}(h_n)}=
 \sum_{j=0}^{n-k+1} \frac{w_{n+2, j+k+1}(h_k)\, t^{j}}{w_{n,j+k-1}(h_n)}. \eqno(10.4)$$
$y_{n,k}(t) $ is a monic polynomial of degree $n+1-k$ and 
	 $$y_{n,n-2}(t)=- \frac{(h_{n-2} - h_{n-1})(h_{n-2} - h_{n+1} )}{(h_n -h_{n-3}) (h_n- h_{n-1})} + \frac{h_{n-2}-h_{n+1}}{h_{n-1}-h_n} t - \frac{h_{n-2}-h_{n+1}}{h_{n-1}-h_n}t^2 +t^3. $$ 
Using (6.3) and (6.4) we obtain
$$y_{n,n-2}(t)= -1 + z t - z t^2 + t^3, \qquad n \ge 2. $$ 
This is the characteristic polynomial of the difference equation (6.1). Let us denote it by $y(t)$.
We will show in the next Lemma that $y(t)$ divides $y_{n,k}(t)$ for $n \ge 2$ and $ 0 \le k \le n-2$. Let us suppose that 
 this result holds.

Let $E$ denote the forward shift operator defined by $E s_k= s_{k+1}$ for any sequence $s_k$ and $k \ge 0.$
Then $y_{n,k}(E)$ is a linear operator that acts on sequences and 
$$y_{n,k}(E)d_k=\sum_{j=k}^{n+1} \frac{w_{n+2,j+1}(h_k)}{w_{n,j-1}(h_n)} E^{j-k} d_k= \hat R_{n,k}. \eqno(10.5)$$
Since the characteristic polynomial $y(t)$ of the difference equation (6.1) divides $y_{n,k}(t)$, for $n\ge 2$ and $0 \le k \le n-2$, $d_k$ satisfies (6.1), and $d_0=0$, it is clear that $\hat R_{n,k}=0$ for $n\ge 0$ and $ 0 \le k \le n-2.$
Therefore $\hat R$ is tridiagonal and hence $R$ is tridiagonal. Note that the condition $d_0=0$ is needed in (10.5) in the case $k=0$. 

Suppose now that $R$ is tridiagonal. Then $\hat R$ is tridiagonal and since $y(t)=y_{n,n-2}(t)$ for $n\ge 2$, from (10.5) we obtain 
$y(E) d_{n-2}=0$ for $n \ge 2$, and this means that the sequence $d_k$ satisfies the difference equation (6.1).
This completes the proof of the theorem. \eop

\begin{lemma} The characteristic polynomial $y(t)$ of the difference equation (6.1) divides the polynomials $y_{n,k}(t)$, defined in (10.4), for $n \ge 2$ and $ 0 \le k \le n-2$.
\end{lemma}

{\it Proof:} 
Define the polynomials 
$$ \nu_{n,k}(t)= \frac{1}{w_{k+2,k-1}(h_n)}\sum_{j=0}^{n-k-2} \frac{w_{n+2,j+k+1}(h_k)\ t^j}{w_{n,j+k+2}(h_n)}, \qquad n \ge 2,\ 0 \le k \le n-2. \eqno(10.6)$$
We claim that $y_{n,k}(t) = - y(t) \nu_{n,k}(t)$ for $n \ge 2,\ 0 \le k \le n-2.$
Let us write the product $ y(t) \nu_{n,k}(t)$ as 
$$ y(t) \nu_{n,k}(t)= \sum_{j=0}^{n-k+1} \rho(n,k,j) \, t^j, \eqno(10.7) $$
and the characteristic polynomial $y(t)$ of (6.1) as 
$$y(t)=\psi_3 t^3 + \psi_2 t^2 + \psi_1 t + \psi_0 = t^3 - z t^2 + z t -1. $$ 
The coefficient of $t^j$ in $ y(t) \nu_{n,k}(t)$ is given by
$$\rho(n,k,j)= \sum_{i=\max(0,j-3)}^{\min(j,n-k-2)} \psi_{j-i} \, \frac{w_{n+2,i+k+1}(h_k)}{ w_{k+2,k-1}(h_n) \, w_{n,i+k+2}(h_n)},\eqno(10.8) $$
We divide $\rho(n,k,j)$ by the coefficient of $t^j$ in $y_{n,k}(t)$ and, after a simplification, we obtain the function 
$$\xi(n,k,j) = \sum_{i=\max(0,j-3)}^{\min(j,n-k-2)} \psi_{j-i} \, \frac{w_{j+k+1,i+k+1}(h_k) w_{i+k+2,j+k-1}(h_n)}{w_{k+2,k-1}(h_n)}. \eqno(10.9)$$
Now we must show that $\xi(n,k,j)=-1$ for $n \ge 2,\ 0 \le k \le n-2, \ 0\le j \le n-k+1.$ Let us note that the common denominator in (10.9) depends only on $n$ and $k$ and it is a polynomial of degree 3 evaluated at $h_n$.
For $j=n-k+1$ we obtain
$$\xi(n,k,n-k+1)= \frac{w_{n+2,n-1}(h_k)}{w_{k+2,k-1}(h_n)}, $$
which equals $-1$ by (6.4). A straightforward computation shows that $\xi(n,k,n-k+1) -\xi(n,k,n-k)$ is a multiple of $ z (h_{n-1}-h_{n}) - h_{n-2} + h_{n+1}$, which is equal to zero, since $h_k$ satisfies (6.1). In the same way we find that $\xi(n,k,n-k) -\xi(n,k,n-k-1)$ is a multiple of $z(h_{n-2}-h_{n-1}) -h_{n-3}+h_{n}$, which is also zero.

We have analogous results at the other end of the values of $j$.
For $j=0$ we obtain $\xi(n,k,0)=-1.$ We define $h_{-1}= h_2 -z (h_1 -h_0)$, to avoid undefined terms in (10.9). 
We find that $ \xi(n,k,1) - 1 $ is a multiple of $z (h_{k} - h_{k+1}) - h_{k-1} + h_{k+2}$, which equals zero, and also that 
$\xi(n,k,2)-\xi(n,k,1) $ is a multiple of $ z (h_{k+1} - h_{k+2}) - h_{k} + h_{k+3}$, which is zero.

It remains to show that $\xi(n,k,j)=-1$ for $3 \le j \le n-k-2$. Such terms exist only when $n-k\ge 5$.
Define the polynomials
$$ \zeta(n,k,j,t)= - \sum_{i=j-3}^{j} \psi_{j-i} \, w_{j+k+1,i+k+1}(h_k) w_{i+k+2,j+k-1}(t). \eqno(10.10)$$
 This is obtained from (10.9) by eliminating the denominator, replacing $h_n$ with $t$ and changing sign. 
 It is convenient to change the summation index. Put $\ell =i-j+3$. Then we have 
$$ \zeta(n,k,j,t)= - \sum_{\ell=0}^{3} \psi_{3-\ell} \, w_{j+k+1,\ell +j +k-2}(h_k) w_{\ell +j +k-1,j+k-1}(t). \eqno(10.11)$$
We will show that $\zeta(n,k,j,t)$ is equal to $w_{k+2,k-1}(t)$ for $n \ge 5,\  0 \le k \le n-5,\ 3 \le j \le n-k-2.$

Let $m \ge 0$. The set $\{ w_{m+i,m}(t): 0 \le i \le 3\}$ is a Newton type basis for the space of polynomials of degree at most 3. Newton's interpolation formula \cite[Eq. 2.10, Exer. 2.2-1]{CdB}, \cite[Eq. 2.26]{ddci}, gives us 
$$ w_{k+2,k-1}(t) = \sum_{\ell=0}^3 \sum_{i=0}^\ell \frac{ w_{k+2,k-1}(h_{m+i})}{w_{m+\ell+1, m}^\prime(h_{m+i })} w_{m+\ell,m}(t). \eqno(10.12)$$
Taking $m= j+k-1$ we obtain 
$$ w_{k+2,k-1}(t) = \sum_{\ell=0}^3 \sum_{i=0}^\ell \frac{ w_{k+2,k-1}(h_{j+k-1 +i})}{w_{j+k+\ell,j+k-1 }^\prime(h_{j+k-1 +i })} w_{\ell +j+k-1,j+k-1}(t). \eqno(10.13)$$
It remains to show that corresponding coefficients of the polynomials $ w_{\ell +j+k-1,j+k-1}(t)$ in (10.11) and (10.13) coincide. The coefficient for $\ell=3$ in (10.11) is equal to 1. The corresponding coefficient in (10.13) is equal to 1 by a general property of divided differences. It is the coefficient of the monic Newton polynomial of degree three in the expression for a monic polynomial of degree three in terms of the Newton polynomials.

For $\ell=0$ the coefficient in (10.11) is $-w_{j+k+1,j+k-2}(h_k)$ and the corresponding coefficient in (10.13) is $w_{k+2,k-1}(h_{j+k-1}) $ and by (6.4) their quotient is equal to 1. 

 The coefficient for $\ell=1$ in (10.13) divided by $ z ( h_k- h_{j+k-1}) (h_k - h_{j+k})$, which is the corresponding coefficient in (10.11), is  
$$ \frac{w_{k+2,k-1}(h_{j+k-1}) - w_{k+2,k-1}(h_{j+k})}{z( h_k- h_{j+k-1}) (h_k - h_{j+k}) (h_{j+k-1}-h_{j+k})} = \qquad \qquad \qquad $$
$$ \frac{1}{z( h_{j+k-1} - h_{j+k})} \left( \frac{-(h_{j+k-1}-h_{k-1}) (h_{j+k-1}-h_{k+1})}{(h_k-h_{j+k})} + \frac{(h_{j+k}-h_{k-1}) (h_{j+k}-h_{k+1})}{(h_{k}-h_{j+k-1})}\right). $$
The first term in the last sum can be written as
$$ \frac{(h_{j+k-1}-h_{k-1}) (h_{j+k-1}-h_{k+1})}{(h_{j+k-2}-h_{k})(h_{j+k}-h_{k})} (h_{j+k-2}-h_{k})= h_{j+k-2}-h_{k},$$
and the second term can be written as  
$$\frac{-(h_{j+k}-h_{k-1}) (h_{j+k}-h_{k+1})}{(h_{j+k-1}- h_k) (h_{j+k+1}-h_k)} (h_{j+k+1}-h_k)=-( h_{j+k+1}-h_k).$$ 
The quotients in the last two equations are equal to 1 and -1 respectively, by (6.4).
Therefore the quotient of the coefficients for $\ell=1$ is equal to 
$$\frac{ h_{j+k-2}-h_{k}- ( h_{j+k-1}-h_k)}{z ( h_{j+k-1} - h_{j+k})} = \frac{ h_{j+k-2}- h_{j+k+1}}{z ( h_{j+k-1} - h_{j+k})} = 1.$$

For $\ell=2$ the corresponding coefficient in (10.11) is equal to $-z (h_k- h_{j+k})$. 
By a general property of divided differences the coefficient for $\ell=2$ in (10.13) is equal to 
$$h_{j+k+1} +h_{j+k}+ h_{j+k-1 } - h_{k+1} - h_k - h_{k-1} =\sum_{i=k-1 }^{j+k-2} (h_{i+3} - h_i).$$
This can also be obtained by simplification of the sum with three terms that corresponds to $\ell =2$ in (10.13).
Since $h_{i+3} - h_i= z(h_{i+2} - h_{i+1} ) $, the sum in the previous equation is equal to the telescopic sum
$$ z \sum_{i=k-1 }^{j+k-2} (h_{i+2} - h_{i+1}) = z( h_{j+k} - h_k) ,$$ 
and therefore it is equal to the corresponding term for $\ell=2$ in (10.11). Therefore we have proved that $y(t)$ divides $y_{n,k}(t)$, for $n\ge 2$ and $0 \le k \le n-2$. \eop

In the third step we use ``diagonal" similarity to show that the matrix $L$ is tridiagonal if and only if $R$ is tridiagonal.

\begin{theorem}
	Let $h_k$, $x_k$, and $d_k$ be sequences that satisfy the difference equation (6.1) and such that $d_0=0$ and $h_j \ne h_k$ if $ j \ne k$.
	Define the sequence $g_k = x_{k-1} (h_k - h_0) + d_k$, for $k \ge 1$, and $g_0=0$, and suppose that $g_k \ne 0$ for $k\ge 1.$
	Then the matrix 
	$$L = C \, (S\tr + F)\, C^{-1} $$ is tridiagonal, where $C$ is the matrix whose entries are 
	$$ c_{n,k} =\prod_{j=k}^{n-1} \frac{ g_{j+1}}{ h_n - h_j}, \qquad 0 \le k \le n-1,$$
	 and $c_{n,n}=1$ for $n\ge 0$, $F$ is the diagonal matrix whose $(k,k)$ entry is $x_k$ for $k \ge 0$, and $S\tr$ is the right shift matrix, defined in (3.5).
\end{theorem}
	
{\it Proof:} Let $G$ be the diagonal matrix whose $(k,k)$ entry is $g_{k+1}$ for $k \ge 0$. We will show first that  
$L_g = K\, ( G \, S\tr + F) \, K^{-1}$ is tridiagonal.

We have 
$$(K \, F \, K^{-1})_{n,k}=\sum_{j=k}^n \frac{x_j}{w_{n,j}(h_n) w_{j+1,k+1}(h_k)}=\sum_{j=k+1}^{n+1} \frac{x_{j-1}}{w_{n,j-1}(h_n) w_{j,k+1}(h_k)}. \eqno(10.14) $$
In the sum
$$\sum_{j=k}^{n+1} \frac{x_{j-1} ( h_j -h_k)}{w_{n,j-1}(h_n) w_{j+1,k+1}(h_k)}, $$
the summand corresponding to $j=k$ is zero, and if $j>k$ the negative of the factor $h_j -h_k$ appears as a factor  
of $w_{j+1,k+1}(h_k)$. Therefore this sum is the negative of the sum in (10.14). Using this result and writing 
 $g_j=x_{j-1} (h_j - h_k + h_k - h_0) + d_j$, for $j \ge 1$, we obtain 
$$(L_g)_{n,k}= (K \, ( G \, S\tr +F)\, K^{-1})_{n,k} = \sum_{j=k}^{n+1} \frac{x_{j-1} ( h_k -h_0) + d_j}{w_{n,j-1}(h_n) w_{j+1,k+1}(h_k)} .\eqno(10.15)$$
Since $x_j$ and $d_j$ satisfy (6.1), $d_0=0$, and $h_k - h_0$ is independent of the summation index $j$, if we put $x_{-1}=x_2-z(x_1-x_0)$ then we can apply Theorem 10.1 and we conclude that $L_g$ is tridiagonal.

Define the $g$-factorial function $f_g$ by $f_g(0)=1$ and $f_g(k)=g_1 g_2 \cdots g_k$ for $k\ge 1.$
Let $\Gamma$ be the diagonal matrix whose $(k,k)$ entry is $f_g(k)$ for $k \ge 0$. Since $g_k \ne 0$, for $k \ge 1$, it is clear that $\Gamma$ is invertible and that the $(k,k)$ entry of $\Gamma^{-1}$ is $1/f_g(k)$, for $ k \ge 0$.

It is easy to verify that 
$$ C= \Gamma \, K \, \Gamma^{-1}, \quad S\tr = \Gamma \, G \, S\tr \, \Gamma^{-1}, \quad F= \Gamma \, F \, \Gamma^{-1}. \eqno(10.16)$$
Therefore 
$$ \Gamma \, L_g \, \Gamma^{-1} = \Gamma \, K \, ( G \, S\tr + F ) \, K^{-1} \, \Gamma^{-1}= C \, (S\tr +F)\, C^{-1} = L. \eqno(10.17)$$ Since $L_g$ is tridiagonal so is $L$. This completes the proof of the theorem. \eop

Combining the results obtained in this Section we obtain a proof of Theorem 6.1.

\section{Final remarks}
If we know the coefficients $c_{n,k}$ with respect to a Newton basis of some orthogonal polynomial sequence, and if we also know the sequence of eigenvalues $\{h_k\}$ then the sequence $\{g_k\}$ can be computed using the recurrence 
$$ c_{n,k}= \frac{g_{k+1}}{h_n - h_k} \  c_{n,k+1}, \qquad 0 \le k \le n-1, $$ 
 which we used to find the explicit formula for $c_{n,k}$.

The number of parameters required to describe particular subfamilies of orthogonal or $q$-orthogonal polynomials can be reduced if we take into account the invariance properties of the coefficients $\alpha_k$ and $\beta_k$ of the three-term recurrence relation, rescaling of the independent variable $t$, the parameters that are zero, and the translations of the polynomials. 

In sections 8 and 9 we have used the normalized recurrence coefficients to determine the parameters for the particular families because they do not depend on the basis of the space of polynomials used to express the polynomial sequences, and also because the procedure to find the parameters for a given family can also be used to see if there are several Newton basis that can be used for the given polynomial family. The parameters $a_0,a_1,a_2,b_0,b_1,b_2,d_1,d_2$ provide a uniform description of all the families and are useful for the purpose of classification.

It seems that the approach used in \cite{Tsuji} to study symmetric polynomials may be used to complement our approach so that we can deal with Hermite and $q$-Hermite polynomials.

The class $\cH_{-1}$ is another source of research problems. It seems that only a few examples of elements in that class have been studied. See \cite{Vinet1}, \cite{Vinet2}, and \cite{Vinet3}.

The classification of the families in the Askey and the $q$-Askey scheme using the uniform parametrization is another interesting research topic.

\section*{Acknowledgements}
The author wishes to thank an anonymous referee that made a meticulous examination of the manuscript and provided numerous ideas to improve the present paper. The author also thanks Professor Tom Koornwinder for his valuable comments about previous versions of this paper.

\end{document}